\documentclass[11pt,a4paper]{article}

\usepackage{latexsym,amsfonts,amstext,amsbsy,amsmath,amssymb,amscd,graphics,epsfig,eucal}

\usepackage{pb-diagram}
\usepackage[all]{xy}
\usepackage{makeidx}
\usepackage{indentfirst}
\usepackage{graphicx}
\usepackage{mathrsfs}

\begin{document}

\date{}
\def\e{\varepsilon}
\def\w{\widetilde}
\def\b{\Box}

\begin{center}
{\large \bf A PACKAGE OF PROGRAMS FOR DETERMINATION OF SOME CLASSES OF SUBGROUPOIDS}\\[0.2cm]
\end{center}
\begin{center}
{\bf Gheorghe Ivan and Mihai Ivan} \\[0.1cm]
\end{center}
\begin{center}
{\it West University of Timi\c soara, Seminarul de Algebr\u a, 4,
Bd. V. P{\^a}rvan, 300223, Timi\c soara, Romania\\
( e-mail: ivan @ math.uvt.ro)}
\end{center}

{\small {\bf Abstract.}  In this paper we give three programs on
computer for finding the subgroupoids, wide subgroupoids and normal
subgroupoids of a finite groupoid. Applying these programs for
groups, we can determine all subgroups and normal subgroups of a
given finite group.
}.\\

\begin{center}
{\bf 1. INTRODUCTION}\\
\end{center}

The concept of groupoid has introduced by H. Brandt [Math. Ann. {\bf
96}, 360 -366 ( 1926; JFM 52.0110.09 )]. A groupoid is an algebraic
structure determined by a partially composition law and a nonempty
set of units. In the language of categories, a groupoid is a small
category in which all morphisms are invertible. For more details and
references about groupoids the reader can be consult the papers (
[1]-[3], [6]-[8]).

The plan of this paper is as follows. In the first section we give
some preliminary concepts concerning groupoids. In the second
section we present an algorithm for finding the subgroupoids of a
groupoid. This algorithm is implemented on computer and we obtain
the program $ BGroidAP2. $  This program has published in [4; Zbl
1109.20310]. In the Section $3$ and Section $4$ we give the programs
$ BGroidAP3 $ and $ BGroidAP4 $  for to determine the wide
subgroupoids resp. the normal subgroupoids of a finite groupoid. We
illustrate the utilization of these programs on some finite
groupoids.

The programs exposed in this paper represent an essential tool for the study of finite groupoids.\\

\begin{center}
{\bf 2. CLASSES OF SUBGROUPOIDS}\\
\end{center}

Let $( G,G_{0})$ be a pair of nonempty sets with $G_{0}\subseteq G,$
endowed with the surjections $ \alpha, \beta : G \to G_{0} $, called
the {\it source}  and the {\it target} map, respectively, a {\it (
partial ) composition  law} $ \mu :G_{(2)}\longrightarrow G, (x,y)
\longrightarrow  \mu(x,y), $ where $ G_{(2)} =\{ (x,y)\in G\times G
| \beta (x)=\alpha (y) \}$ and an injection $\iota : G \to G, x \to
\iota(x)$ called the {\it inversion map}. We write sometimes $
x\cdot y $ or $ x y $ for $ \mu(x,y) $ and $ x^{-1} $ for $
\iota(x).$ The elements of $G_{(2)}$ are called {\it composable
pairs} of $ G. $

{\it DEFINITION 2.1.} ([4]) $(i)$ The $5$-tuple  $( G, \alpha,
\beta, \mu ; G_{0} )$ is a {\it semigroupoid}, if the composition
law is {\it associative}, i.e. $ (x y) z = x (y z), $ for all $ x,
y, z \in G $ such that the products $ (xy)z $ and $ x(yz) $ are
defined.

$(ii) $ A {\it monoidoid} is a semigroupoid $ (G, \alpha, \beta,
\mu ; G_{0}) $ such that the {\it identities property} holds, i.e.
for each $ x\in G $ we have $ (\alpha (x),x), (x,\beta (x))\in
G_{(2)} $ and $ \alpha(x) x = x\beta(x) = x.$

$(iii) $ The $6$- tuple $ ( G, \alpha, \beta, \mu, \iota; G_{0} )
$ is a {\it groupoid} or
 a {\it $ G_{0}$-groupoid}, if\\
$ (G,\alpha, \beta, \mu ; G_{0}) $ is a monoidoid such that
 the {\it inverses property} holds, i.e. for each $ x\in G $ we have $ (x^{-1},x), (x,x^{-1})\in G_{(2)} $
 and $ x^{-1} x = \beta(x),  x x^{-1} = \alpha(x). $ \hfill$\b$

{\it  Remark 2.1.} The definition of the groupoid is equivalent as
the one used in the paper ([2]).\hfill$\b$

The element $ \alpha(x) $ [resp. $ \beta(x) $ ] denoted sometimes
by $ u_{l}(x) $ [ resp. $ u_{r}(x) $] is the {\it left unit} [
resp. {\it right unit}]  of $ x\in G.$ The set $ G_{0} $ is called
the {\it unit set} of $ G. $ A $ G_{0} $-groupoid $ G $ will be
denoted by $ ( G,\alpha,\beta; G_{0}) $ or $ ( G ; G_{0}). $ The
maps $ \alpha, \beta, \mu  $ and $ \iota $ are called the {\it
structure functions} of $ G. $

If $(G,\alpha,\beta; G_{0})$ is a groupoid, then the following properties hold (see [3]):\\
$(1)   \alpha(u) = \beta(u) = u  ,  u\cdot u = u $  and $ \iota(u)=u $ for all $ u\in G_{0}; $\\
$(2)   \alpha(x y) = \alpha(x) $ and $  \beta(x y ) = \beta(y),  (\forall)  (x,y)\in G_{(2)}; $\\
$(3)    \alpha (x^{-1}) = \beta (x),    \beta ( x^{-1}) = \alpha (x) $  for all $ x\in G; $\\
$(4)   G(u)=\{ x\in G | \alpha(x)=\beta(x)=u \}$ is a group under
the restriction of $\mu $ to $ G(u),$ called the {\it isotropy
group at} $ u $ of $ G.$\hfill$\b$

{\it Example 2.1.} $(i)$ A group $ G $ having $ e $ as unity, is
just a $ \{ e\} $- groupoid in the following way: the maps $ \alpha,
\beta : G \to G_{0} $ and $ \iota : G \to G $ are
 given by $ \alpha(x)=\beta(x)=e , \iota(x)= x^{-1} $ for all $ x\in G; $
for all $ x,y\in G $ the element $ x\cdot y $ is the product of
elements $ x $ and $ y $ in the group $ G.$ Conversely, every
groupoid $ G $ with one unit is a group.

$(ii) $ The {\it nul groupoid over a set}. Any nonempty set $ X $
may be regarded as a groupoid on itself with the groupoid
structure : $ G = G_{0} = X, \alpha = \beta =\iota= Id_{X}; $ $ x
, y\in X $ are composable iff $  x=y  $ and we define $  x\cdot x
= x.  $\hfill$\b$

{\it Example 2.2.} $(i)$ {\it The groupoid $ {\cal F}_{inj}(S,X) $}.
For a nonempty set $ X $ denote by ${\cal F}_{inj}(S,X) =\{ f : S
\to X |  (\forall) S, \emptyset\neq S\subseteq X, f \hbox{is
injective} \}.$ For $ f\in {\cal F}_{inj}(S,X), $ let $ D(f) $ be
the domain of $ f $ and let $ R(f) = f(D(f)). $  For $ G={\cal
F}_{inj}(S,X), $ let $ G_{(2)}=\{(f,g)\in G\times G | R(f)=D(g) \} $
and for $ (f,g)\in G_{(2)} $ define $ \mu(f,g)=g\circ f. $ If $
Id_{S} $ denotes the identity map on $ S, $ then $ G_{0} = \{ Id_{S}
| \emptyset \neq S\subseteq X \} $ is the set of units of $ G. $ The
maps $ \alpha, \beta : G \to G_{0} $ and $ \iota : G \to G $ are
defined by $ \alpha(f) = Id_{D(f)}, \beta(f)=Id_{R(f)} $ and $
\iota(f)= f^{-1}. $ Thus $ {\cal F}_{inj}(S,X) $ is a groupoid,
called the {\it groupoid of injective functions from the nonempty
sets $ S $ of $ X $ into $ X $}.

In particular, if $ X = \{ 1, 2,\ldots, n \}, $  the groupoid of
the injective functions defined on the subsets of $ \{ 1,2,\ldots,
n \} $ is called the {\it symmetric groupoid of degree $ n $} and
is denoted by $ {\cal S}_{n}; $ for several properties of $ {\cal
S}_{n}, $ see [5].

$(ii) $  Let $ Oxy $ be a system of cartesian coordinates in a
plane. We consider the subsets $ Ox = \{(x,0)\in {\bf R}^{2} |
(\forall) x\in {\bf R} \} $ and $ Oy= \{(0,y)\in {\bf R}^{2} |
(\forall) y\in {\bf R} \} $ of $ X={\bf R}^{2}. $ Let $ G =\{
f_{1}=Id_{Ox}, f_{2}=Id_{Oy}, f_{3}=\sigma_{Ox} ,
f_{4}=\sigma_{Oy} \} \subset {\cal F}_{inj}(S,{\bf R}^{2}) $ where
$ f_{3} : Ox \to Oy, f_{4} : Oy \to Ox $ are defined by $
f_{3}(x,0)=(0,x) $ and $ f_{4}(0,y)=(y,0) $ ($ \sigma_{Ox} $ resp.
$ \sigma_{Oy} $ is called the {\it saltus function defined on} $ x
$ {\it -axis } resp. $ y $-{\it axis} ).

For the composable pairs of $ G $, the map $ \mu $ is defined by:
$\mu(f_{1},f_{1})=f_{1};$ $\mu(f_{1},f_{3})= f_{3}\circ
f_{1}=f_{3};$ $\mu(f_{2},f_{2})= f_{2};$  $\mu(f_{2},f_{4})=
f_{4};$ $\mu(f_{3},f_{2})= f_{3};$ $\mu(f_{3},f_{4})=f_{1};$
$\mu(f_{4},f_{1})=f_{4};$ $\mu(f_{4},f_{3})=f_{2}. $

The unit set of $ G $ is $ G_{0} =\{ f_{1}=Id_{Ox},f_{2}=Id_{Oy}
\} $ and  $ \alpha, \beta : G \to G_{0} , \iota : G \to G $ are
given by $ \alpha(f_{j})=\beta(f_{j})=\iota(f_{j})=f_{j} $ for $
j=1,2 ; \alpha(f_{3})=\beta(f_{4})=f_{1};
 \alpha(f_{4})=\beta(f_{3})=f_{2};  \iota(f_{3})=f_{4} $ and
$ \iota(f_{4})=f_{3}. $  It is easy to verify that $ (G; G_{0}) $
is a groupoid, denoted by $ {\cal F}_{(4;2)}({\bf R}^{2}) $ and
called the {\it groupoid of saltus functions defined on the axes
of coordinates in a plane}.\hfill$\b$

{\it Example 2.3.} If $\{G_{i} | i\in I\}$ is a disjoint family of
groupoids, let $ G = \cup_{i \in I} G_{i},
 G_{(2)}=\cup_{i \in I} G_{i,(2)}$ and $ G_{0} = \cup_{i \in I} G_{i,0}$, where $ G_{i,0}$ is the unit set of $ G_{i}$.
 Here, $ x,y \in G$ may be composed iff they lie in the same groupoid
$ G_{i}$ and they are composable in $ G_{i}$. This groupoid is
denoted by $ \coprod_{i\in I} G_{i} $ and is called the {\it
disjoint union of groupoids} $ G_{i}, i\in I. $ In particular, the
disjoint union of groups $ G_{i}, i \in I$ is a
groupoid.\hfill$\b$

A finite $ G_{0}-$ groupoid $ G $ such that $ | G | = n $ and $ |
G_{0} | = m $ is called {\it $ ( n ; m )-$ groupoid} or {\it
finite groupoid of type $ (n;m).$} We will sometimes denote a
finite groupoid of type $(n;m)$ by $G_{(n;m)}.$

{\it Example 2.4.} $(i) $ Each finite groupoid of type $ (n;1) $ is
a group.

$(ii) $ Each finite groupoid of type $ (n;n) $ is a nul groupoid.

$(iii) $ The groupoid $ {\cal F}_{(4;2)}({\bf R}^{2}) $  is a $
(4;2)-$ groupoid.\hfill$\b$

{\it DEFINITION 2.2.} $ (i) $ Let $ (G, \alpha, \beta;G_{0}) $ be a
groupoid. A pair $ (H;H_{0}) $ of nonempty sets such that $ H
\subseteq G $ and $ H_{0}\subseteq G_{0} $ is a {\it subgroupoid} of
$ G, $ if
 the following conditions hold: $ (1)   \alpha(H)=\beta(H)=H_{0} ; $
$ (2) $ for all $ x,y\in H $ such that $ x y $ is defined, we have
$ x y\in H $ and $ (3) $ for all $ x\in H, $ we have $ x^{-1}\in
H.$

$(ii) $ A subgroupoid $ (H;H_{0}) $ of a $ G_{0}-$ groupoid $ G $
with property that $ H_{0} = G_{0} $ is called {\it wide
subgroupoid} of $ G. $

$(iii) $ A wide subgroupoid $ (H; G_{0}) $ of a groupoid $ (
G;G_{0}) $
 is called a {\it normal subgroupoid} of $ G, $ if for all $ x\in G $ and $ h\in H $ such that
 the product $ x h x^{-1} $ is defined, we have $ x h x^{-1} \in H. $\hfill$\b$

The intersection of any collection of subgroupoids of a groupoid is itself a subgroupoid
 of that groupoid. If $ G $ is a $ G_{0}-$ groupoid and $ X $ is a nonempty subset of $ G, $ then
 the intersection of all subgroupoids of $ G $ which contain $ X $ is a subgroupoid, denoted
 by $< X >$ and called the {\it generated subgroupoid} of $ G $ by $ X.$

{\it Example 2.5.} $ (i) $ In a group $ G, $ every subgroupoid (in
fact, subgroup) is a wide subgroupoid and conversely.

$(ii) $ If $ G $ is a $ G_{0} $- groupoid, then $ G_{0} $ is a
normal subgroupoid of $ G, $ called the {\it nul subgroupoid} of $
G.$

$(iii) $ Let be the Klein $ 4 $-group $ K_{4} = \{ (1), \sigma =
(12)(34), \tau=(13)(24), \sigma\circ \tau = (14)(23) \}\subset
S_{4} $  ( it is a subgroup of the symmetric group $ S_{4} $ of
degree $ 4 $). We have $ \sigma^{2} = \tau^{2} = (1) $ and $
\tau\circ \sigma =\sigma\circ \tau.$ We consider the disjoint
union $ G = K_{4}\coprod {\cal F}_{(4;2)}({\bf R}^{2}) $ of the
group $ K_{4} $ with the groupoid $ {\cal F}_{(4;2)}({\bf R}^{2})
$ given in Example 1.2 (ii). We have that $G=\{ (1), \sigma ,
\tau, \sigma\circ \tau,
f_{1}=Id_{Ox},f_{2}=Id_{Oy},f_{3}=\sigma_{Ox},f_{4}=\sigma_{Oy} \}
$ is a $ G_{0}$- groupoid of type $(8;3),$ where $ G_{0} = \{ (1),
f_{1},f_{2} \}.$  It is easy to verify that $K_{4} $ and $ {\cal
F}_{(4;2)}({\bf R}^{2})$  are subgroupoids of $ G.$ Also, $ N_{1}
= K_{4}\coprod \{ f_{1}, f_{2} \} $ and $ N_{2} = \{(1)\}\coprod
{\cal F}_{(4;2)}({\bf R}^{2}) $
 are wide subgroupoids of $ G.$ Moreover, $ N_{1} $ and $ N_{2} $ are normal subgroupoids of $ G. $\hfill$\b$\\

\begin{center}
{\bf  $3.$ ALGORITHM FOR DETERMINATION OF SUBGROUPOIDS. THE $ BGroidAP 2 $  PROGRAM}\\
\end{center}

We consider a given finite universal algebra $ (G, \alpha, \beta,
\mu , \iota;G_{0}) $ such that $ | G | = n $ and $ | G_{0} | =
m $ with $ 1\leq m\leq n. $ We denote the elements of $ G $ by\\
$ a_{1}, a_{2}, \cdots, a_{m}, a_{m+1}, \cdots, a_{n} $ such that
$ G_{0} =\{ a_{1}, a_{2}, \cdots, a_{m} \}.$

We give an algorithm for decide if the universal algebra $
(G,\alpha, \beta, \mu, \iota; G_{0}) $ is a $ G_{0} $- groupoid
and for determine the subgroupoids of $ G. $ This algorithm is
constituted by the following stages.

{\bf Stage I.} {\it We introduce the initial data}: $ n =| G |, m
= | G_{0} |;$ the functions $ \alpha,
 \beta, \iota $ and $ \mu $ given by its tables of structure.

{\bf Stage II.} {\it Test if the universal algebra $ (G, \alpha,
\beta, \mu, \iota; G_{0}) $  is a groupoid}. For this, the following steps are executed:\\
{\bf step 1.} $ (G,\alpha, \beta, \mu,\iota ; G_{0} ) $ is a
structure well-defined, i.e. $ \alpha, \beta $ are surjections,
$ \iota $ is injective and $ \mu $ is defined on $ G_{(2)} $ with values in $ G; $\\
{\bf step 2.} $(G,\alpha, \beta, \mu ; G_{0}) $ is a semigroupoid;\\
{\bf step 3.}  the semigroupoid $ (G, \alpha, \beta, \mu; G_{0}) $ is a monoidoid;\\
{\bf step 4.} the monoidoid $ (G,\alpha,\beta,\mu; G_{0} ) $ is a groupoid.\\
{\bf step 5.} If the above steps are satisfied, make the tables of
the structure functions $ \alpha, \beta, \iota $ and $ \mu $ and
write the message "{\it $ G $ is a groupoid}$\,$".

{\bf Stage III.} {\it Determine the subgroupoids of} $ G. $
The following steps must be executed:\\
{\bf step 1.}$   $ Write all nonempty subsets $ X $ of $ G $;\\
{\bf step 2.}$   $ Determine the subgroupoid $ < X > $ of $ G $ generated by $ X; $\\
{\bf step 3.}$   $ Sort by cardinal all subgroupoids determined in the step $ 2;$\\
{\bf step 4.}$   $ List the subgroupoids produced in the above step;\\
{\bf step 5.}$   $ For each subgroupoid  make its subgroupoid table.\\[-0.2cm]

Let us we present the correspondence between the initial data and
input data:\\[-0.2cm]
$$\begin{array}{ccc}
G=\{ a_{1},a_{2},\ldots,\ldots,a_{m},a_{m+1},\ldots,a_{n} \} &
\longleftrightarrow & \{ 1,2,\ldots,m,m+1,\ldots,n \}\cr Initial~
data &\longleftrightarrow & Input~  data\cr | G |=n
&\longleftrightarrow & n\cr | G_{0} |=m &\longleftrightarrow & m\cr
\end{array}$$
$$\begin{array}{|c|c|c|c|c|c|c|}\hline
 a_{k}    & \!a_{1}\!\!\! & \!\!\cdots\!\! & \!a_{m}\!\! & a_{m+1} & \!\!\cdots\!\! & a_{n} \\ \hline
  \!\alpha(a_{k})\!\!& \!a_{1}\!\!\! & \!\!\cdots\!\! & \!a_{m}\!\! & \alpha(a_{m+1}) & \!\!\cdots\!\! &\!\alpha(a_{n})\!\! \cr \hline
 \! \beta(a_{k})\!\!&\! a_{1}\!\!\! &\!\! \cdots\!\! & \!a_{m}\!\! & \beta(a_{m+1}) & \!\!\cdots\!\! & \beta(a_{n}) \cr \hline
\!\iota(a_{k})\!\! & \!a_{1}\!\!\! & \!\!\cdots\!\! & \!a_{m}\!\!& \iota(a_{m+1}) & \!\!\cdots \!\!& \iota(a_{n}) \cr \hline
\end{array} \longleftrightarrow \begin{array}{cccccc} \\
       1 \!\!&\!\! \cdots\!\! & m & \!\!u_{l}(m+1)\!\! & \!\!\cdots\!\! & \!\!u_{l}(n)\cr
       1 \!\!&\!\! \cdots\!\! & m & \!\!u_{r}(m+1)\!\! & \!\!\cdots\!\! & \!\!u_{r}(n) \cr
       1 \!\!&\!\! \cdots\!\! & m & \!\!inv(m+1) \!\!& \!\!\cdots\!\! & \!\!inv(n)\cr
\end{array}$$
$$\begin{array}{|r|c|c|c|c|c|}\hline
\mu  & a_{1} & \cdots & a_{k}            & \cdots & a_{n} \\\hline
a_{1}  &       &        &                  &        &\cr \hline
\cdots &       &        &                  &        &\cr \hline
a_{j}  &       &        & a_{j}\cdot a_{k} &        &\cr \hline
\cdots &       &        &                  &        &\cr \hline
a_{n}  &       &        &                  &        &\cr \hline
\end{array}  \longleftrightarrow  \begin{array}{cccccc} \\
a_{11}  & \cdots    &  a_{1k}   & \cdots  & a_{1n} \cr
a_{21}  & \cdots    & a_{2k}    & \cdots  &  a_{2n}\cr
\cdots  & \cdots    &\cdots    & \cdots  & \cdots  \cr
a_{j1}  &  \cdots   & a_{jk} & \cdots  & a_{jn} \cr
\cdots  & \cdots    & \cdots    & \cdots  &\cdots  \cr
a_{n1}  & \cdots    & a_{nk}    & \cdots  & a_{nn}.\cr
\end{array}$$

The absence of an element from the arrow $ "j" $ and the column $
"k" $ of the table of $ \mu $ indicates the fact that the pair $
(a_{j}, a_{k})\in G\times G $ is not composable. The element $
a_{jk}=\mu(a_{j},a_{k}) $ is represented by $ 0 $ in the table of
input data, if the product $ a_{j}\cdot a_{k} $ is not defined.

{\it Example 3.1.} Let the groupoid $ G = K_{4}\coprod {\cal
F}_{(4;2)}({\bf R}^{2}), $ see Example 2.5 (iii). We have $ G=\{
a_{1} = (1), a_{2} = f_{1} , a_{3}= f_{2}, a_{4} = \sigma, a_{5} =
\tau, a_{6} = \sigma\circ \tau , a_{7} = f_{3}, a_{8} = f_{4}\} $
and the correspondence between the initial data and
input data are the following:\\[-0.2cm]
$$\begin{array}{ccc}
G=\{ a_{1},a_{2},a_{3},a_{4},a_{5}, a_{6}, a_{7}, a_{8} \} &
\longleftrightarrow & \{ 1,2, 3, 4, 5, 6, 7 , 8 \}\cr Initial~ data
&\longleftrightarrow & Input~  data\cr | G |=8 &\longleftrightarrow
& 8\cr | G_{0} |=3 &\longleftrightarrow & 3\cr
\end{array}$$\\[-0.5cm]
$$\begin{array}{|c|c|c|c|c|c|c|c|c|}\hline
 a_{k}        & a_{1} & a_{2} & a_{3} & a_{4} & a_{5} & a_{6} & a_{7} & a_{8}\\ \hline
\alpha(a_{k}) & a_{1} & a_{2} & a_{3} & a_{1} & a_{1} & a_{1} &  a_{2}    & a_{3} \cr \hline
\beta(a_{k})  & a_{1} & a_{2} & a_{3} & a_{1} & a_{1} & a_{1} &  a_{3}  & a_{2} \cr \hline
\iota(a_{k})  & a_{1} & a_{2} & a_{3} & a_{4} & a_{5} & a_{6} & a_{8} & a_{7} \cr \hline
\end{array} \longleftrightarrow \begin{array}{ccccccccc} \\
       1  & 2 & 3 & 1 & 1 & 1 & 2 & 3\cr
       1  & 2 & 3 & 1 & 1 & 1 & 3 & 2 \cr
       1  & 2 & 3 & 4 & 5 & 6 & 8 & 7 \cr
\end{array}$$\\[-0.5cm]
$$\begin{array}{|r|c|c|c|c|c|c|c|c|}\hline
\mu    & a_{1} & a_{2} & a_{3} & a_{4} & a_{5} & a_{6} & a_{7} & a_{8}\\ \hline
a_{1}  & a_{1} &       &       & a_{4} & a_{5} & a_{6} &       &       \cr \hline
a_{2}  &       & a_{2} &       &       &       &       & a_{7} &       \cr \hline
a_{3}  &       &       & a_{3} &       &       &       &       & a_{8}  \cr \hline
a_{4}  & a_{4} &       &       & a_{1} & a_{6} & a_{5} &       &  \cr \hline
a_{5}  & a_{5} &       &       & a_{6} & a_{1} & a_{4} &       &  \cr \hline
a_{6}  & a_{6} &       &       & a_{5} & a_{4} & a_{1} &       &  \cr \hline
a_{7}  &       &       & a_{7} &       &       &       &       & a_{2} \cr \hline
a_{8}  &       & a_{8} &       &       &       &       & a_{3} &     \cr \hline
\end{array}  \longleftrightarrow  \begin{array}{cccccccc} \\
1 & 0 & 0 & 4 & 5 & 6 & 0 & 0\cr
0 & 2 & 0 & 0 & 0 & 0 & 7 & 0\cr
0 & 0 & 3 & 0 & 0 & 0 & 0 & 8\cr
4 & 0 & 0 & 1 & 6 & 5 & 0 & 0\cr
5 & 0 & 0 & 6 & 1 & 4 & 0 & 0\cr
6 & 0 & 0 & 5 & 4 & 1 & 0 & 0\cr
0 & 0 & 7 & 0 & 0 & 0 & 0 & 2\cr
0 & 8 & 0 & 0 & 0 & 0 & 3 & 0\cr
\end{array}$$

The implementation of the above algorithm on computer is realized
in the program $ BGroidAP2, $ which is composed from two modules
denoted by $ unit21.dfm $ and $ unit21.pas $. The module $
unit21.pas $ is consists from the principal program followed of
procedures and functions.

The principal program of the module $ unit21.pas $ is constituted from the following lignes.\\

\begin{tabular}{|c|c|l|}\hline
& Lignes & The module $ unit21.pas $ \\ \hline\hline & 001 &  unit
Unit1;\cr & 002 & interface\cr & 003 & uses\cr & 004 &
\hspace*{0.5cm}Windows, Messages, SysUtils, Classes, Graphics,
Controls,\cr & &\hspace*{0.5cm}Forms, Dialogs, Grids, DBGrids,
ShellAPI, Db, DBTables,\cr &  & \hspace*{0.5cm}StdCtrls, Menus,
ExtCtrls, ComCtrls, ToolWin, Spin;\cr & 005 & const\cr & 006
&\hspace*{1.3cm}nmax = 200;\cr & 007 & type\cr & 008
&\hspace*{1.3cm}TSubSet = Set of Byte;\cr & 009
&\hspace*{0.5cm}TForm1 = class(TForm)\cr & 010
&\hspace*{0.8cm}MainMenu1: TMainMenu;\cr & 011
&\hspace*{0.8cm}File1: TMenuItem;\cr & 012
&\hspace*{0.8cm}OpenFile1: TMenuItem;\cr & 013
&\hspace*{0.8cm}SaveFile1: TMenuItem;\cr & 014
&\hspace*{0.8cm}GroupBox1: TGroupBox;\cr & 015
&\hspace*{0.8cm}StringGrid1: TStringGrid; \cr & 016
&\hspace*{0.8cm}StringGrid2: TStringGrid;\cr & 017
&\hspace*{0.8cm}GroupBox2: TGroupBox;\cr & 018
&\hspace*{0.8cm}StringGrid3: TStringGrid; \cr & 019
&\hspace*{0.8cm}StringGrid4: TStringGrid;\cr & 020
&\hspace*{0.8cm}OpenDialog1: TOpenDialog;\cr & 021
&\hspace*{0.8cm}SaveDialog1: TSaveDialog;\cr & 022
&\hspace*{0.8cm}Splitter1: TSplitter; \cr & 023
&\hspace*{0.8cm}Splitter2: TSplitter;\cr & 024
&\hspace*{0.8cm}ToolBar1: TToolBar;\cr & 025
&\hspace*{0.8cm}ToolBar2: TToolBar;\cr & 026
&\hspace*{0.8cm}Splitter4: TSplitter; \cr & 027
&\hspace*{0.8cm}ToolButton1: TToolButton; \cr & 028
&\hspace*{0.8cm}ToolButton2: TToolButton; \cr & 029
&\hspace*{0.8cm}New1: TMenuItem;\cr & 030
&\hspace*{0.8cm}ToolBar3: TToolBar;\cr & 031
&\hspace*{0.8cm}ToolButton4: TToolButton;\cr & 032
&\hspace*{0.8cm}ToolButton5: TToolButton;\cr & 033
&\hspace*{0.8cm}Label2: TLabel;\cr \hline
\end{tabular}

\begin{tabular}{|c|c|l|}\hline
& Lignes & The module $ unit21.pas $ \\ \hline\hline & 034
&\hspace*{0.8cm}SpinEdit1: TSpinEdit;\cr & 035
&\hspace*{0.8cm}ToolButton6: TToolButton;\cr & 036
&\hspace*{0.8cm}Label3: TLabel;\cr & 037
&\hspace*{0.8cm}SpinEdit2: TSpinEdit;\cr & 038
&\hspace*{0.8cm}Savesubgroupoid1: TMenuItem;\cr & 039
&\hspace*{0.8cm}StatusBar1: TStatusBar;\cr & 040
&\hspace*{0.8cm}StatusBar2: TStatusBar;\cr & 041
&\hspace*{0.8cm}ToolButton3: TToolButton;\cr & 042
&\hspace*{0.8cm}ToolButton11: TToolButton;\cr & 043
&\hspace*{0.8cm}procedure FormShow(Sender: TObject);\cr & 044
&\hspace*{0.8cm}procedure Button1Click(Sender: TObject);\cr
$\star$ & 045 &\hspace*{0.8cm}procedure Button2Click(Sender:
TObject);\cr $\star$ & 046 &\hspace*{0.8cm}procedure
StringGrid1SetEditText(Sender: TObject;\cr & & \hspace*{1.1cm}
ACol,ARow: Integer; const Value: String);\cr $\star$ & 047
&\hspace*{0.8cm}procedure StringGrid2SetEditText(Sender:
TObject;\cr & &\hspace*{1.1cm}ACol, ARow: Integer; const Value:
String);\cr $\star$ & 048 &\hspace*{0.8cm}procedure
OpenFile1Click(Sender: TObject);\cr $\star$ & 049
&\hspace*{0.8cm}procedure SaveFile1Click(Sender: TObject);\cr &
050 &\hspace*{0.8cm}procedure StringGrid3SelectCell(Sender:
TObject; ACol,\cr &&\hspace*{1.1cm}ARow: Integer; var CanSelect:
Boolean);\cr $\star$ & 051 &\hspace*{0.8cm}procedure
New1Click(Sender: TObject);\cr $\star$& 052
&\hspace*{0.8cm}procedure ToolButton4Click(Sender: TObject);\cr &
053 &\hspace*{0.8cm}procedure Savesubgroupoid1Click(Sender:
TObject);\cr $\star$ & 054 &\hspace*{0.8cm}procedure
ToolButton3Click(Sender: TObject);\cr & 055
&\hspace*{0.5cm}private \cr & 056 &\hspace*{1.4cm}ForcedStop :
Boolean;\cr & 057 &\hspace*{0.8cm}err\_message : String;\cr & 058
&\hspace*{1.4cm}\qquad subgr : array[1..10000] of TSubSet;\cr &
059 &\hspace*{0.8cm}units, SelectedSub : TSubSet;\cr & 060
&\hspace*{1.4cm}m, n, nsub : Integer;\cr & 061 &\hspace*{0.8cm}h :
array[0..nmax, 0..nmax] of Byte;\cr & 062 &\hspace*{0.8cm}u\_left,
u\_right, inv : array[0..nmax] of Integer;\cr $\star$ & 063
&\hspace*{0.8cm}procedure WMDropFiles(var Msg: TWMDropFiles);\cr
&&\hspace*{1.1cm} message WM\_DROPFILES;\cr $\star$ & 064
&\hspace*{0.8cm}procedure PerformFileOpen(const FileName1 :
string);\cr $\star$ & 065 &\hspace*{0.8cm}procedure
PerformFileSave(const FileName1 : string);\cr
  & 066 &\hspace*{0.8cm}procedure PerformSaveSubgroupoid(t : TSubSet;\cr
&&\hspace*{1.1cm} const FileName1 : string);\cr
$\star$ & 067 &\hspace*{0.8cm}procedure MakeUnitsTable;\cr
$\star$ & 068 &\hspace*{1.4cm}\qquad procedure MakeGroupoidTable;\cr
 & 069 &\hspace*{0.8cm}procedure MakeSubgroupoidTable(t : TSubSet);\cr
$\star$ & 070 &\hspace*{0.8cm}function ToStr(x : Integer) : String;\cr
 \hline
\end{tabular}

\begin{tabular}{|c|c|l|}\hline
& Lignes & The module $ unit21.pas $ \\ \hline\hline & 071
&\hspace*{0.8cm}function SubsetToString(t : TSubSet) : String;\cr
 & 072 &\hspace*{0.8cm}function Cardinal(t : TSubSet) : Byte;\cr
 & 073 &\hspace*{0.8cm}procedure Cover(var t : TSubSet);\cr
 & 074 &\hspace*{0.8cm}function AlreadyFound(t : TSubSet) : Boolean;\cr
 & 075 &\hspace*{0.8cm}procedure AddSubgroupoid(t : TSubSet);\cr
 & 076 &\hspace*{0.8cm}procedure GenerateSubgroupoids(t : TSubSet; r : Byte);\cr
 & 077 &\hspace*{0.8cm}procedure SortByCardinal;\cr
 & 078 &\hspace*{0.8cm}procedure ListSubgroupoids;\cr
$\star$ & 079 &\hspace*{0.8cm}function IsStructure : Boolean;\cr
$\star$ & 080 &\hspace*{0.8cm}function IsSemigroupoid : Boolean;\cr
$\star$ & 081 &\hspace*{0.8cm}function IsMonoidoid : Boolean;\cr
$\star$ & 082 &\hspace*{0.8cm}function IsGroupoid : Boolean;\cr
& 083 &\hspace*{0.5cm}public\cr
& 084 &\hspace*{0.5cm}end;\cr
& 085 &var\cr
& 086 &\hspace*{0.5cm}Form1: TForm1;\cr
& 087 & implementation\cr
& 088 &\{\$R *.DFM\}\cr
& 089 &procedure TForm1.FormShow(Sender: TObject);\cr
& 090 &var\cr
& 091 &\hspace*{1.4cm}i, j : Byte;\cr
& 092 &begin\cr
& 093 &\hspace*{0.5cm}DragAcceptFiles(Handle, True);\cr
& 094 & \hspace*{0.5cm}StringGrid1.EditorMode := True;\cr
& 095 & \hspace*{0.5cm}n := 0;\cr
& 096 & \hspace*{0.5cm}m := 0;\cr
& 097 & \hspace*{0.5cm}for i := 0 to nmax do\cr
& 098 & \hspace*{1.4cm}for j := 0 to nmax do\cr
& 099 & \hspace*{1.4cm}h[i,j]:= 0\cr
& 100 &\hspace*{0.5cm}for i := 0 to nmax do begin\cr
& 101 &\hspace*{0.8cm}u\_left[i] := 0;\cr
& 102 &\hspace*{0.8cm}u\_right[i] := 0;\cr
& 103 &\hspace*{0.8cm}inv[i] := 0;\cr
& 104 &\hspace*{0.5cm}end;\cr
& 105 &\hspace*{0.5cm}nsub := 0;\cr
& 106 &\hspace*{0.5cm}SelectedSub := [ ]\cr
& 107 & end;\cr
& 108 & end.\cr\hline
\end{tabular}\\

The procedures and functions marked by the symbol $" \star "$ can be
find in [4] or in the preprint arXiv:math/0602604v1 [math GR]. The
other procedures and
functions contained in the module $ unit21.pas $ are presented in the follows.\\[-0.2cm]

{\it procedure TForm1.Cover;}

var

\hspace*{1.4cm}i, j : Byte;

\hspace*{0.5cm}modif : Boolean;

begin

\hspace*{0.5cm}repeat

\hspace*{3.1cm}modif := false;

\hspace*{0.8cm}for i := 1 to n do if i in t then begin

\hspace*{1cm}if not (inv[i] in t) then begin

\hspace*{1.4cm}modif := true;

\hspace*{1.4cm}t := t + [inv[i]]

\hspace*{1cm}end;

\hspace*{1cm}for j := 1 to n do if j in t then

\hspace*{1.4cm}if u\_right[i] = u\_left[j] then

\hspace*{3.1cm}if not (h[i, j] in t) then begin

\hspace*{3.1cm}modif := true;

\hspace*{3.1cm}t := t + [h[i, j]]

\hspace*{1.7cm}end

\hspace*{0.8cm}end

\hspace*{0.5cm}until not modif

end;\\[-0.3cm]

{\it function TForm1.Cardinal;}

var

\hspace*{1.4cm}i, nr : Byte;

begin

\hspace*{1.4cm}nr := 0;

\hspace*{0.5cm}for i := 1 to n do

\hspace*{1.4cm}if i in t then

\hspace*{1.4cm}nr := nr + 1;

\hspace*{0.5cm}Cardinal := nr;

end;\\[-0.3cm]

{\it function TForm1.SubsetToString;}

var

\hspace*{1.4cm}s : String;

\hspace*{0.5cm}i : Byte;

begin

\hspace*{1.4cm}s := '\{';

\hspace*{0.5cm}for i := 1 to n do

\hspace*{1.4cm}if i in t then

\hspace*{1.4cm}s := s + \{'a' +\} tostr(i) + ', ';

\hspace*{0.5cm}delete(s, length(s)-1, 2);

\hspace*{0.5cm}s := s + '\}';

\hspace*{0.5cm}SubsetToString := s;

end;\\[-0.3cm]

{\it procedure TForm1.AddSubgroupoid;}

begin

\hspace*{1.4cm}nsub := nsub + 1;

\hspace*{0.5cm}subgr[nsub] := t;

\hspace*{0.5cm}if nsub = 3000 then

\hspace*{1.4cm}ForcedStop := true;

end;\\[-0.3cm]

{\it function TForm1.AlreadyFound;}

var

\hspace*{1.4cm}i : Integer;

begin

\hspace*{1.4cm}AlreadyFound := False;

\hspace*{0.5cm}for i := 1 to nsub do

\hspace*{1.4cm}if t = subgr[i] then

\hspace*{1.4cm}AlreadyFound := True;

\hspace*{0.5cm}if t = [ ] then

\hspace*{1.4cm}AlreadyFound := True

end;\\[-0.3cm]

{\it procedure TForm1.GenerateSubgroupoids;}

var

\hspace*{1.4cm}i : Byte;

begin

\hspace*{1.4cm}Cover(t);

\hspace*{1.4cm}if not AlreadyFound(t) then

\hspace*{1.4cm}AddSubgroupoid(t);

\hspace*{0.5cm}for i := r to n do

\hspace*{1.4cm}if not (i in t) then

\hspace*{1.4cm}if not ForcedStop then

\hspace*{1.4cm}GenerateSubgroupoids(t + [i], i);

end;\\[-0.3cm]

{\it procedure TForm1.SortByCardinal;}

var

\hspace*{1.4cm}i, j : Integer;

\hspace*{0.5cm}aux : TSubSet;

begin

\hspace*{1.4cm}for i := 1 to nsub - 1 do

\hspace*{1.4cm}for j := i + 1 to nsub do

\hspace*{1.4cm}if Cardinal(subgr[i]) $ > $ Cardinal(subgr[j]) then begin

\hspace*{1.4cm}aux := subgr[i];

\hspace*{1.4cm}subgr[i] := subgr[j];

\hspace*{1.4cm}subgr[j] := aux

\hspace*{1cm}end

end;\\[-0.3cm]

{\it procedure TForm1.ListSubgroupoids;}

var

\hspace*{1.4cm}i : Integer;

begin

\hspace*{0.5cm}StringGrid3.RowCount := nsub;

\hspace*{0.5cm}for i := 1 to nsub do

\hspace*{1.7cm}StringGrid3.Cells[0, i - 1] := SubsetToString(subgr[i])

end;\\[-0.3cm]

{\it procedure TForm1.PerformSaveSubgroupoid;}

var

\hspace*{1.4cm}f : TextFile;

\hspace*{0.5cm}i, j, ordin, nunits : Byte;

\hspace*{0.5cm}sir, ind : array[1..nmax] of Byte;

begin

\hspace*{0.5cm}ordin := 0;

\hspace*{0.5cm}for i := 1 to n do

\hspace*{1.4cm}if i in SelectedSub then begin

\hspace*{1.4cm}ordin := ordin + 1;

\hspace*{1cm}sir[ordin] := i;

\hspace*{1cm}ind[i] := ordin

\hspace*{0.8cm}end;

\hspace*{0.5cm}nunits := cardinal(SelectedSub * units);

\hspace*{1.4cm}AssignFile(f, FileName1);

\hspace*{0.5cm}rewrite(f);

\hspace*{0.5cm}writeln(f, ordin);

\hspace*{0.5cm}writeln(f, nunits);

\hspace*{0.5cm}writeln(f);

\hspace*{0.5cm}for i := 1 to ordin do

\hspace*{1.4cm}write(f, ind[u\_left[sir[i]]], ' ');

\hspace*{0.5cm}writeln(f);

\hspace*{0.5cm}for i := 1 to ordin do

\hspace*{1.4cm}write(f, ind[u\_right[sir[i]]], ' ');

\hspace*{0.5cm}writeln(f);

\hspace*{0.5cm}for i := 1 to ordin do

\hspace*{1.4cm}write(f, ind[inv[sir[i]]], ' ');

\hspace*{0.5cm}writeln(f);

\hspace*{0.5cm}writeln(f);

\hspace*{0.5cm}for i := 1 to ordin do begin

\hspace*{1.4cm}for j := 1 to ordin do

\hspace*{1.4cm}write(f, ind[h[sir[i], sir[j]]], ' ');

\hspace*{0.8cm}writeln(f);

\hspace*{0.5cm}end;

\hspace*{0.5cm}CloseFile(f)

end;\\[-0.3cm]

{\it procedure TForm1.MakeSubgroupoidTable;}

var

\hspace*{0.5cm}ordin, i, j : Byte;

\hspace*{1.4cm}sir : array[1..nmax] of Byte;

begin

\hspace*{0.5cm}ordin := 0;

\hspace*{0.5cm}for i := 1 to n do

\hspace*{1.4cm}if i in t then begin

\hspace*{1.4cm}ordin := ordin + 1;

\hspace*{1cm}sir[ordin] := i

\hspace*{0.8cm}end;

\hspace*{0.5cm}StringGrid4.RowCount := ordin + 1;

\hspace*{0.5cm}StringGrid4.ColCount := ordin + 1;

\hspace*{0.5cm}for i := 1 to ordin do begin

\hspace*{1.4cm}StringGrid4.Cells[0, i] := tostr(sir[i]);

\hspace*{1.4cm}StringGrid4.Cells[i, 0] := tostr(sir[i])

\hspace*{0.5cm}end;

\hspace*{0.5cm}for i := 1 to ordin do

\hspace*{1.4cm}for j := 1 to ordin do

\hspace*{1.4cm}if h[sir[i], sir[j]] $ < > $ 0 then

\hspace*{3.1cm}StringGrid4.Cells[j, i] := tostr(h[sir[i], sir[j]])

\hspace*{1cm}else

\hspace*{3.1cm}StringGrid4.Cells[j, i] := ''

end;\\[-0.3cm]

{\it procedure TForm1.Savesubgroupoid1Click(Sender: TObject);}

begin

\hspace*{1.4cm}if SelectedSub $ < > $ [ ] then

\hspace*{3.1cm}if SaveDialog1.Execute then

\hspace*{3.1cm}PerformSaveSubgroupoid(SelectedSub, SaveDialog1.FileName)

\hspace*{0.8cm}else

\hspace*{0.5cm}else Application.MessageBox('No subgroupoid selected.', '', mb\_OK)

end;\\[-0.2cm]

{\it procedure TForm1.StringGrid3SelectCell(Sender: TObject; ACol,

\hspace*{0.5cm}ARow: Integer; var CanSelect: Boolean);}

begin

\hspace*{0.5cm}MakeSubgroupoidTable(subgr[ARow + 1]);

\hspace*{0.5cm}SelectedSub := subgr[ARow + 1]

end;\\[-0.3cm]

{\it procedure TForm1.Button1Click(Sender: TObject);}

begin

\hspace*{1.4cm}nsub := 0;

\hspace*{1.4cm}ForcedStop := false;

\hspace*{0.5cm}GenerateSubgroupoids([ ], 1);

\hspace*{0.5cm}SortByCardinal;

\hspace*{1.4cm}ListSubgroupoids;

\hspace*{1.4cm}StatusBar2.SimpleText := tostr(nsub) + ' subgroupoid(s) found.'

end;.\\[-0.3cm]

We illustrate the utilization of the program $ BGroidAP2 $ in the
following cases.

{\it Example 3.2.} {\it Determination of subgroupoids of a groupoid
of type $ (8;2). $} We consider the subset $ G_{(8;2)}=\{ g_{j} |
j=\overline{1,8} \}$ of the symmetric groupoid ${\cal S}_{3} $ ( see
[5] ), where: $g_{1}=\left (\begin{array}{cc}
1 & 2\\
1 & 2
\end{array}\right),$
$g_{2}=\left (\begin{array}{cc}
1 & 3\\
1 & 3
\end{array}\right),$
$g_{3}=\left (\begin{array}{cc}
1 & 2\\
2 & 1
\end{array}\right),$
$g_{4}=\left (\begin{array}{cc}
1 & 2\\
1 & 3
\end{array}\right),$
$g_{5}=\left (\begin{array}{cc}
1 & 2\\
3 & 1
\end{array}\right),$
$g_{6}=\left (\begin{array}{cc}
1 & 3\\
1 & 2
\end{array}\right),$
$g_{7}=\left (\begin{array}{cc}
1 & 3\\
2 & 1
\end{array}\right),$
$g_{8}=\left (\begin{array}{cc}
1 & 3\\
3 & 1
\end{array}\right).$

We denote the restrictions of the structure functions $ \alpha,
\beta, \iota $ and the composition law defined on the groupoid
${\cal S}_{3} $ to $ G_{(8;2)} $ by the same symbols. Using the correspondence \\[-0.2cm]
$$\begin{array}{ccc}
G_{(8;2)}=\{ g_{1},g_{2},g_{3}, g_{4},g_{5},g_{6}, g_{7},g_{8} \} & \longleftrightarrow & \{ 1,2,3,4,5,6,7,8 \}\cr
\end{array}$$\\[-0.5cm]
the input data are given by the following tables:
$$\begin{array}{cccccccc} \\
8 & & & & & & &\\
2 & & & & & & &\\[0.1cm]

1 & 2 & 1 & 1 & 1 & 2 & 2 & 2\\
1 & 2 & 1 & 2 & 2 & 1 & 1 & 2 \\
1 & 2 & 3 & 6 & 7 & 4 & 5 & 8\\[0.1cm]
\end{array}     \begin{array}{cccccccc} \\
 1 & 0 & 3 & 4 & 5 & 0 & 0 & 0\\
 0 & 2 & 0 & 0 & 0 & 6 & 7 & 8 \\
 3 & 0 & 1 & 5 & 4 & 0 & 0 & 0\\
 0 & 4 & 0 & 0 & 0 & 1 & 3 & 5 \\
 0 & 5 & 0 & 0 & 0 & 3 & 1 & 4\\
 6 & 0 & 7 & 2 & 8 & 0 & 0 & 0\\
 7 & 0 & 6 & 8 & 2 & 0 & 0 & 0\\
 8 & 0 & 0 & 0 & 0 & 7 & 6 & 2\\
\end{array}$$

Execute the program $ BGroidAP2 $ for $ G_{(8;2)} $ and the window
program of obtained results is presented in the Figure $1$.

\begin{figure}[!htb]
  \centering
  \includegraphics[scale=1]{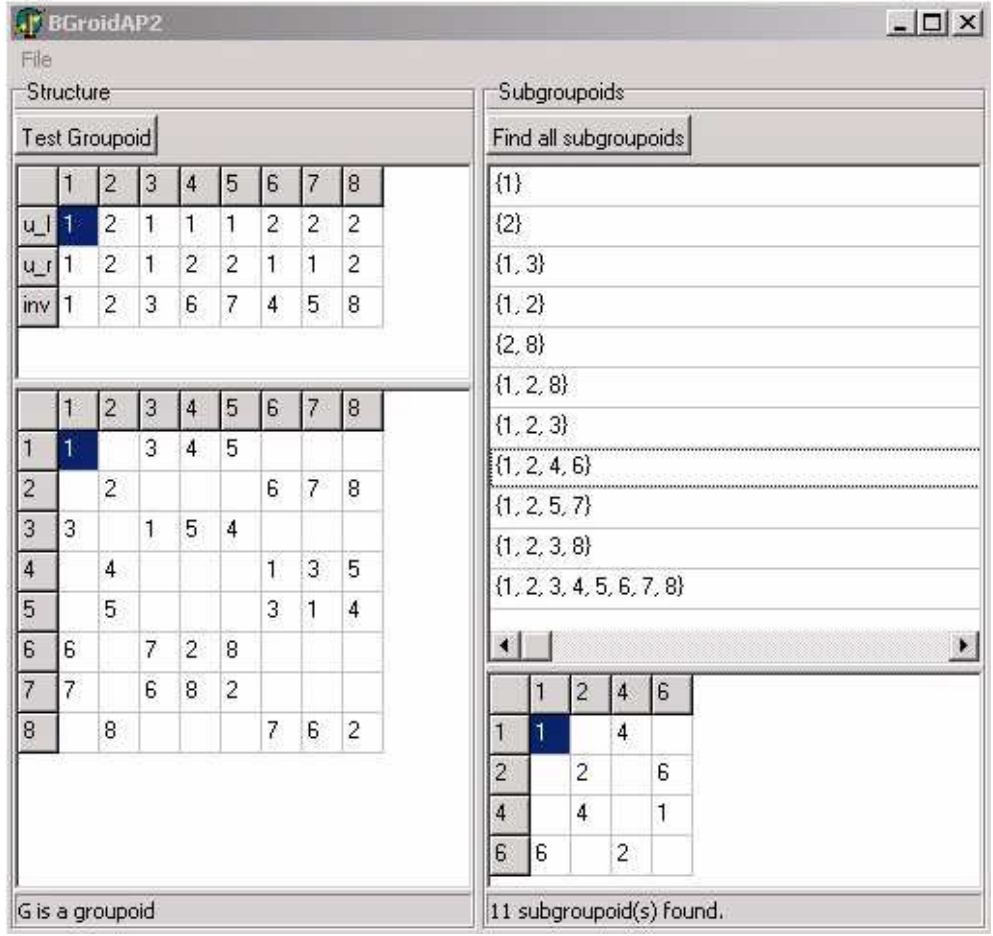}
 \caption{The subgroupoids of a $(8;2)-$ groupoid} \label{fig1}
\end{figure}

Therefore, $ G_{(8;2)} $ is a groupoid with unit set $ G_{(8;2),0} =\{ g_{1}, g_{2}\} $ and it
has $ 11 $ subgroupoids ( see, Fig. 1 ). Using the initial data and input data for this groupoid,
the correspondence between output data and final data is the following :\\[-0.4cm]
$$\begin{array}{lcl}
Output~ data             &\longleftrightarrow & Final~data\cr \{1\},
\{2\}              &\longleftrightarrow & H_{(1;1)}^{1}=\{ g_{1} \},
H_{(1;1)}^{2}=\{ g_{2} \}\cr \{ 1,2\} &\longleftrightarrow &
H_{(2;2)}^{3}=\{ g_{1},g_{2} \}=G_{(8;2),0}\cr \{ 1,3 \}, \{ 2,8 \}
&\longleftrightarrow & H_{(2;1)}^{4}=\{ g_{1},g_{3} \},
H_{(2;1)}^{5}=\{ g_{2},g_{8} \} \cr \{ 1,2,3 \}, \{ 1,2,8 \}
&\longleftrightarrow & H_{(3;2)}^{6}=\{ g_{1}, g_{2}, g_{3} \},
H_{(3;2)}^{7}=\{ g_{1}, g_{2}, g_{8} \}  \cr \{ 1,2,3,8 \}
&\longleftrightarrow & H_{(4;2)}^{8}=\{ g_{1}, g_{2}, g_{3}, g_{8}
\}\cr \{ 1,2,4,6 \} &\longleftrightarrow & H_{(4;2)}^{9}=\{ g_{1},
g_{2}, g_{4}, g_{6} \}\cr \{ 1,2,5,7 \} &\longleftrightarrow &
H_{(4;2)}^{10}=\{ g_{1}, g_{2}, g_{5}, g_{7} \}\cr \{
1,2,3,4,5,6,7,8 \} &\longleftrightarrow & H_{(8;2)}^{11}=
G_{(8;2)}\cr
\end{array}$$\hfill$. \b$
\newpage
\begin{center}
{\bf $ 4.$ DETERMINATION OF WIDE SUBGROUPOIDS. THE $ BGroidAP 3$  PROGRAM}\\
\end{center}

We give an algorithm for decide if the universal algebra $
(G,\alpha, \beta, \mu, \iota; G_{0}) $ is a $ G_{0} $- groupoid
and for determine the all wide subgroupoids of $G. $ This algoritm
is constituted by the following stages.

{\bf Stage I.} {\it We introduce the initial data}, see Section 2.

{\bf Stage II.} {\it Test if the universal algebra $(G, \alpha,
\beta, \mu, \iota; G_{0})$  is a groupoid}. This stage is composed
by five steps; see, Section 2.

{\bf Stage III.} {\it Determine the all wide subgroupoids of} $ G.
$
The following steps must be executed:\\
{\bf step 1.}$   $ Write all nonempty subsets $ X $ of $ G $ with property that $ G_{0}\subseteq X $;\\
{\bf step 2.}$   $ Determine the subgroupoid $ < X > $ of $ G $ generated by $ X; $\\
{\bf step 3.}$   $ Sort by cardinal all wide  subgroupoids determined in the step $ 2;$\\
{\bf step 4.}$   $ List the wide subgroupoids produced in the above step;\\
{\bf step 5.}$   $ For each wide subgroupoid  make its subgroupoid
table.

This algorihtm is implemented on computer and we obtain the program
$ BGroidAP3, $ which is composed from two modules denoted by $
unit31.dfm $ and $ unit31.pas. $

The principal program of the module $ unit31.pas $ consists from the following lignes.\\[0.1cm]

\begin{tabular}{|c|l|}\hline
Lignes & The module $ unit31.pas $ \\ \hline\hline 001 - 027 & the
lignes 001 - 027 of the module unit21.pas;\cr 028 - 036 &  the
lignes 029 - 037 of the module unit21.pas;\cr 037
&\hspace*{0.8cm}ToolButton7: TToolButton;\cr 038 - 043 & the
lignes 038 - 043 of the module unit21.pas;\cr 044 - 051 & the
lignes 045 - 052 of the module unit21.pas;\cr 052
&\hspace*{0.8cm}procedure Button7Click(Sender: TObject);\cr 053 -
107 & the lignes 053 - 107 of the module unit21.pas;\cr 108 &
end.\cr\hline
\end{tabular}\\

The new procedure of the module $ unit31.pas $ is presented in the follows.\\[-0.2cm]

{\it procedure TForm1.ToolButton7Click(Sender: TObject);}

begin

\hspace*{1.4cm}nsub := 0;

\hspace*{1.4cm}ForcedStop := false;

\hspace*{1.4cm}GenerateSubgroupoids(units, m + 1);

\hspace*{0.5cm}SortByCardinal;

\hspace*{1.4cm}ListSubgroupoids;

\hspace*{1.4cm}StatusBar2.SimpleText := tostr(nsub) + ' subgroupoid(s) found.'

end;.\\[-0.3cm]

We illustrate the utilization of the program $ BGroidAP3 $ in the
following cases.

{\it Example 4.1.} {\it Determination of wide subgroupoids of a
disjoint union of two groupoids}. Let the groupoid $ G =
K_{4}\coprod {\cal F}_{(4;2)}({\bf R}^{2}) $ and use the inputs data
presented in the Example 3.1.

Execute the program $ BGroidAP3 $ for $ G $ and the window program
of obtained results is presented in the Figure $2$.

\begin{figure}[!htb]
  \centering
 \includegraphics[scale=1]{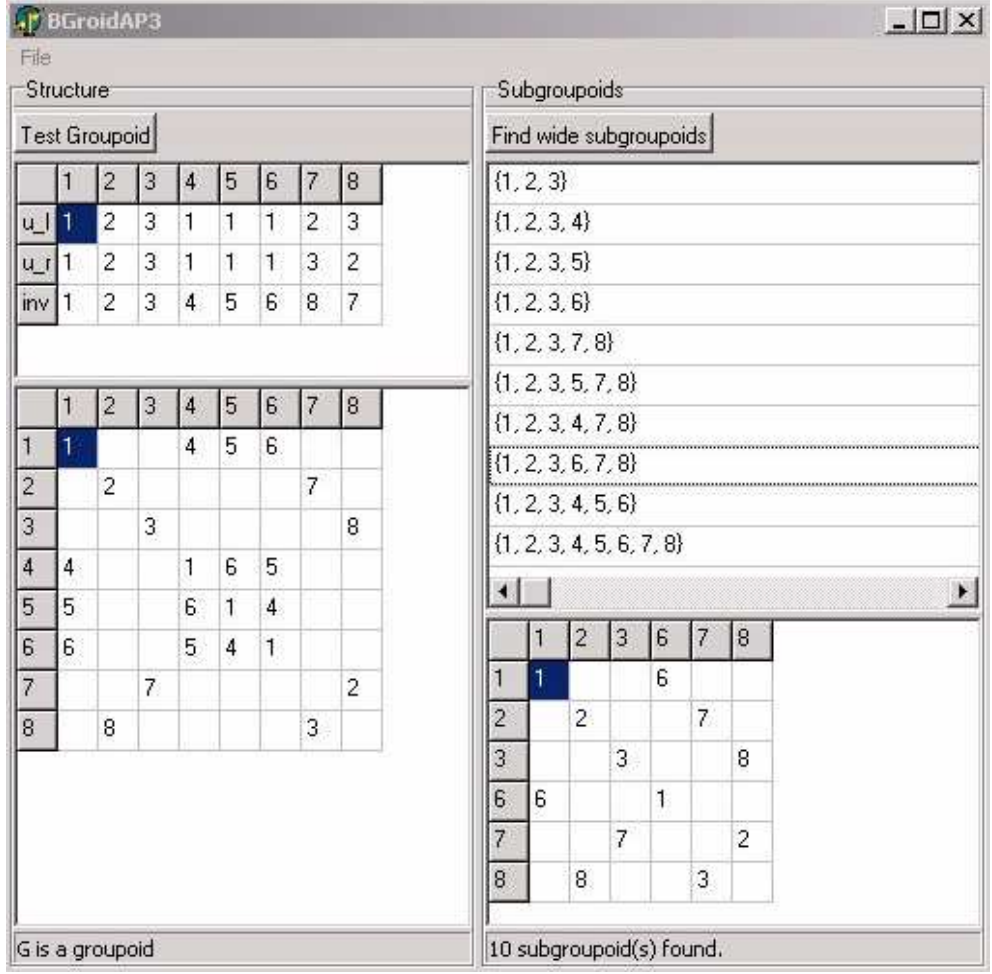}
  \caption{The wide subgroupoids of a $(8;3)-$ groupoid} \label{fig2}
\end{figure}

Therefore, this groupoid has $ 10 $ wide subgroupoids ( see, Fig. 2
). Using the correspondence between output data and initial data,
its wide subgroupoids are the following : $ W_{(3;3)}^{1}=G_{0},$
$W_{(4;3)}^{2}=\{(1),Id_{Ox},Id_{Oy},\sigma\} ,$
$W_{(4;3)}^{3}=\{(1),Id_{0x}, Id_{Oy},\tau\},$
$W_{(4;3)}^{4}=\{(1),Id_{Ox}, Id_{Oy},\sigma\circ\tau\},$
$W_{(5;3)}^{5}=\{ (1), Id_{Ox}, Id_{0y}, \sigma_{Ox},\sigma_{Oy}\},$
$W_{(6;3)}^{6}=\{(1),Id_{Ox},Id_{Oy},\sigma,\tau,\sigma\circ\tau\},$
 $W_{(6;3)}^{7}=\{
 (1),Id_{Ox},Id_{Oy},\sigma,\sigma_{Ox},\sigma_{Oy}\},$
 $W_{(6;3)}^{8}=\{(1),
 Id_{Ox},Id_{Oy},\tau,\sigma_{Ox},\sigma_{Oy}\},$
$W_{(6;3)}^{9}$
$=\{(1),Id_{Ox},Id_{Oy},\sigma\circ\tau,\sigma_{Ox},\sigma_{Oy}\},$
$W_{(8;3)}^{10}= K_{4}\coprod {\cal F}_{(4;2)}({\bf R}^{2}).$

Applying program $ BGroidAP2,$  we obtain that $ G $
 has $ 29 $ subgroupoids.

{\it Example 4.2.} {\it Determination of wide subgroupoids of a
groupoid of type $ (8;2). $} Let the groupoid $ G_{(8;2)}=\{ g_{j} |
j=\overline{1,8} \} $ given in  Example 3.2. Use the inputs data
presented in the Example 2.2 and execute the program $ BGroidAP3. $
We find that this groupoid has $ 7 $ wide subgroupoids, namely: $
W_{(2;2)}^{1}=G_{(8;2),0},$ $W_{(3;2)}^{2}=\{ g_{1}, g_{2}, g_{3}
\},$ $W_{(3;2)}^{3}=\{ g_{1}, g_{2}, g_{8} \},$ $W_{(4;2)}^{4}=\{
g_{1}, g_{2},$ $g_{3}, g_{8} \},$ $W_{(4;2)}^{5}=\{ g_{1}, g_{2},
g_{4}, g_{6} \},$ $W_{(4;2)}^{6}=\{ g_{1}, g_{2}, g_{5}, g_{7} \},$
$W_{(8;2)}^{7}= G_{(8;2)}$\hfill$.\b$\\

{\bf $ 5.$ DETERMINATION OF NORMAL SUBGROUPOIDS. THE $ BGroidAP 4$
PROGRAM}\\[-0.2cm]

We give an algorithm for decide if the universal algebra $
(G,\alpha, \beta, \mu, \iota; G_{0}) $ is a $ G_{0} $- groupoid
and for determine the normal subgroupoids of $ G. $ This algoritm
is constituted by the following stages.

{\bf Stage I.} {\it We introduce the initial data},see Section $ 2
$.

{\bf Stage II.} {\it Test if the universal algebra $ (G, \alpha,
\beta, \mu, \iota; G_{0}) $  is a groupoid}. This stage is
composed by five steps; see, Section 2.

{\bf Stage III.} {\it Determine the normal subgroupoids of} $ G. $
The following steps must be executed:\\
{\bf step 1.}$   $ Write all nonempty subsets $ X $ of $ G $ with property that $ G_{0}\subseteq X $;\\
{\bf step 2.}$   $ Determine the normal subgroupoid $ < X > $ of $ G $ generated by $ X; $\\
{\bf step 3.}$   $ Sort by cardinal all normal subgroupoids determined in the step $ 2;$\\
{\bf step 4.}$   $ List the normal subgroupoids produced in the above step;\\
{\bf step 5.}$   $ For each normal subgroupoid  make its
subgroupoid table.

The program $ BGroidAP4 $ is composed from two modules denoted by
$ unit41.dfm $ and $ unit41.pas. $

The principal program of the module $ unit41.pas $ consists from the following lignes.\\

\begin{tabular}{|c|l|}\hline
Lignes & The module $ unit41.pas $ \\ \hline\hline 001 - 027 & the
lignes 001 - 027 of the module unit21.pas;\cr 028 - 036 &  the
lignes 029 - 037 of the module unit21.pas;\cr 037
&\hspace*{0.8cm}ToolButton8: TToolButton;\cr 038 - 043 & the
lignes 038 - 043 of the module unit21.pas;\cr 044 - 051 & the
lignes 045 - 052 of the module unit21.pas;\cr 052
&\hspace*{0.8cm}procedure Button8Click(Sender: TObject);\cr 053 -
075 & the lignes 053 - 075 of the module unit21.pas;\cr 076
&\hspace*{0.8cm}procedure GenerateNormal(t : TSubSet; r :
Byte);\cr 077 - 082  & the lignes 077 - 082 of the module
unit21.pas;\cr 083 & \hspace*{0.8cm}function IsNormal(t : TSubSet)
: Boolean;\cr 084 - 108  & the lignes 083 - 107 of the module
unit21.pas;\cr 109 & end.\cr\hline
\end{tabular}\\

The new procedures and the function "IsNormal" of the module $ unit41.pas $ are presented in the follows.\\[-0.3cm]

{\it procedure TForm1.GenerateNormal;}

var

\hspace*{1.4cm}i : Byte;

begin

\hspace*{1.4cm}Cover(t);

\hspace*{0.5cm}if IsNormal(t) then

\hspace*{3.1cm}if not AlreadyFound(t) then

\hspace*{3.1cm}AddSubgroupoid(t);

\hspace*{0.5cm}for i := r to n do

\hspace*{1.4cm}if not (i in t) then

\hspace*{1.4cm}if not ForcedStop then

\hspace*{1.4cm}GenerateNormal(t + [i], i);

end;\\[-0.3cm]

{\it function TForm1.IsNormal;}

var

\hspace*{1.4cm}i, j : Byte;

begin

\hspace*{1.4cm}IsNormal := true;

\hspace*{1.4cm}for i := 1 to n do if i in t then

\hspace*{1.4cm}for j := 1 to n do

\hspace*{1.4cm}if (u\_right[j] = u\_left[i]) and (u\_right[j] = u\_right[i]) then

\hspace*{1.4cm}if not (h[h[j, i], inv[j]] in t) then begin

\hspace*{3.1cm}IsNormal := false;

\hspace*{1.7cm}exit;

\hspace*{1.4cm}end;

end;\\[-0.3cm]

{\it procedure TForm1.ToolButton8Click(Sender: TObject);}

begin

\hspace*{1.4cm}nsub := 0;

\hspace*{1.4cm}ForcedStop := false;

\hspace*{1.4cm}GenerateNormal(units, m + 1);

\hspace*{0.5cm}SortByCardinal;

\hspace*{1.4cm}ListSubgroupoids;

\hspace*{1.4cm}StatusBar2.SimpleText := tostr(nsub) + ' subgroupoid(s) found.'

end;.\\[-0.3cm]

We illustrate the utilization of the program $ BGroidAP4 $ in the
following cases.

{\it Example 5.1.} $(i)$ {\it Determination of normal subgroupoids
of a groupoid of type $ (9;3).$} We consider the subset $
K_{(9;3)}=\{ \varphi_{j} | j=\overline{1,9} \}$ of the symmetric
groupoid ${\cal S}_{3}, $  where:
$\varphi_{1}=\left(\begin{array}{c}
1\\
1
\end{array}\right), \varphi_{2}=\left(\begin{array}{c}
2\\
2
\end{array}\right), \varphi_{3}=\left (\begin {array}{c}
3\\
3
\end{array}\right), \varphi_{4}=\left (\begin{array}{cc}
1\\
2
\end{array}\right),\varphi_{5}=\left (\begin{array}{c}
1\\
3
\end{array}\right), \varphi_{6}=\left (\begin{array}{c}
2 \\
1
\end{array}\right), \varphi_{7}=\left (\begin{array}{c}
2 \\
3
\end{array}\right), \varphi_{8}=\left (\begin{array}{c}
3\\
1
\end{array}\right), \varphi_{9}=\left (\begin{array}{c}
3\\
2
\end{array}\right). $

We denote the restrictions of the structure functions $ \alpha, \beta, \iota $ and the composition law defined on the groupoid
$ {\cal S}_{3} $ to $ K_{(9;3)} $ by the same symbols. Using the correspondence\\[-0.3cm]
$$\begin{array}{ccc}
K_{(9;3)}=\{ \varphi_{1},\varphi_{2},\varphi_{3}, \varphi_{4},\varphi_{5},\varphi_{6}, \varphi_{7},\varphi_{8}, \varphi_{9} \} & \longleftrightarrow & \{ 1,2,3,4,5,6,7,8,9 \}\cr
\end{array}$$\\[-0.4cm]
the input data are given by the following tables:\\[-0.4cm]
$$\begin{array}{ccccccccc} \\
9 & & & & & & & &\\
3 & & & & & & & &\\[0.1cm]
1 & 2 & 3 & 1 & 1 & 2 & 2 & 3 & 3\\
1 & 2 & 3 & 2 & 3 & 1 & 3 & 1 & 2\\
1 & 2 & 3 & 6 & 8 & 4 & 9 & 5 & 7\\[0.1cm]
\end{array}     \begin{array}{ccccccccc} \\
 1 & 0 & 0 & 4 & 5 & 0 & 0 & 0 & 0\\
 0 & 2 & 0 & 0 & 0 & 6 & 7 & 0 & 0\\
 0 & 0 & 3 & 0 & 0 & 0 & 0 & 8 & 9\\
 0 & 4 & 0 & 0 & 0 & 1 & 5 & 0 & 0\\
 0 & 0 & 5 & 0 & 0 & 0 & 0 & 1 & 4\\
 6 & 0 & 0 & 2 & 7 & 0 & 0 & 0 & 0\\
 0 & 0 & 7 & 0 & 0 & 0 & 0 & 6 & 2\\
 8 & 0 & 0 & 9 & 3 & 0 & 0 & 0 & 0\\
 0 & 9 & 0 & 0 & 0 & 8 & 3 & 0 & 0\\
\end{array}$$\\[-0.5cm]

Execute the program $ BGroidAP4 $ for $ K_{(9;3)} $ and the window
program of obtained results is presented in the Figure $3$.

\begin{figure}[!htb]
\centering
\includegraphics[scale=1]{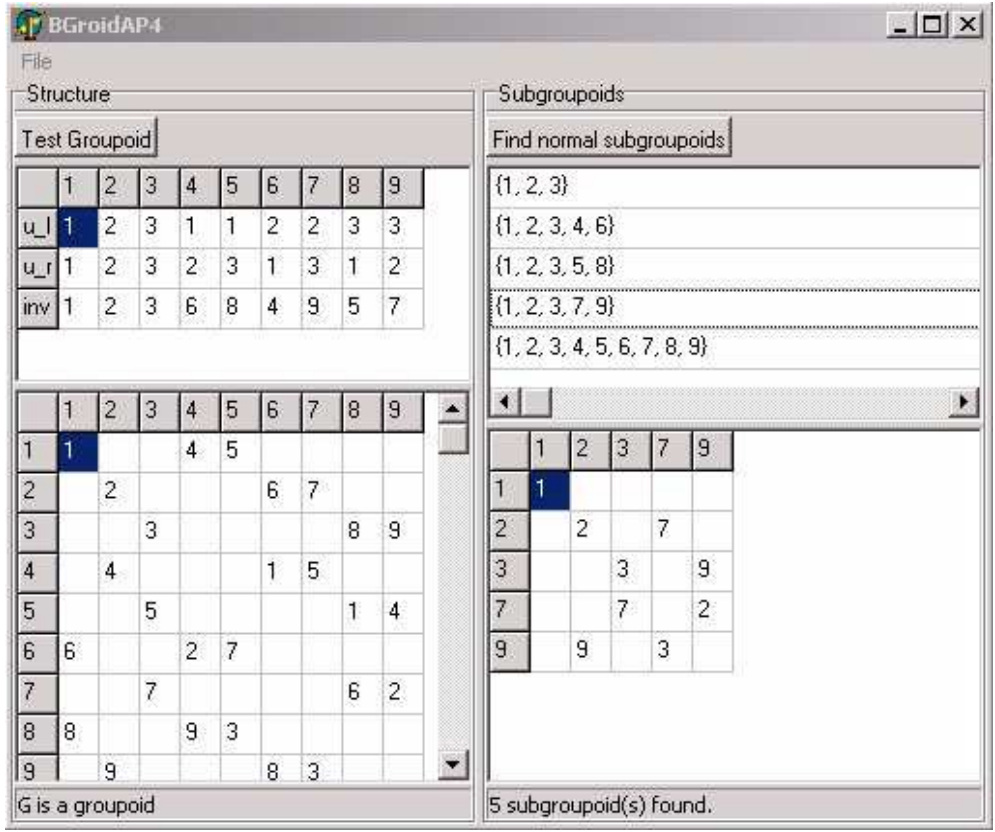}
\caption{The normal subgroupoids of a $(9;3)-$ groupoid}
\label{fig3}
\end{figure}

Therefore, $ K_{(9;3)} $ is a groupoid with unit set $
K_{(9;3),0}=\{ \varphi_{1}, \varphi_{2}, \varphi_{3}\}. $ This
groupoid has $ 5 $ normal subgroupoids ( see, Fig.3 ). Using the
correspondence between output data and initial data, these normal
subgroupoids  are the following : $ N_{(3;3)}^{1}=K_{(9;3),0},$
$N_{(5;3)}^{2}=\{ \varphi_{1}, \varphi_{2}, \varphi_{3},
\varphi_{4},\varphi_{6} \},$ $N_{(5;3)}^{3}=\{ \varphi_{1},
\varphi_{2}, \varphi_{3}, \varphi_{5},\varphi_{8} \},$
$N_{(5;3)}^{4}=\{ \varphi_{1}, \varphi_{2}, \varphi_{3},
\varphi_{7},\varphi_{9} \},$ $N_{(9;3)}^{5}= K_{(9;3)}.$

$(ii) $ Applying the program $ BGroidAP2 $ for the groupoid $
K_{(9;3)} $ we obtain that it has $ 14 $ subgroupoids, namely:
$H_{(1;1)}^{1}=\{\varphi_{1}\},$ $H_{(1;1)}^{2}=\{\varphi_{2}\},$
$H_{(1;1)}^{3}=\{\varphi_{3}\},$
$H_{(2;2)}^{4}=\{\varphi_{1},\varphi_{2}\},$
$H_{(2;2)}^{5}=\{\varphi_{1},\varphi_{3}\},$
$H_{(2;2)}^{6}=\{\varphi_{2},\varphi_{3}\},$
 $H_{(3;3)}^{7}=K_{(9;3),0},$
$H_{(4;2)}^{8}=\{\varphi_{1},\varphi_{2},\varphi_{4},
\varphi_{6}\},$
$H_{(4;2)}^{9}=\{\varphi_{1},\varphi_{3},\varphi_{5},
\varphi_{8}\},$ $H_{(4;2)}^{10}=\{\varphi_{2},\varphi_{3},$
$\varphi_{7}, \varphi_{9}\},$  $H_{(5;3)}^{11}=N_{(5;3)}^{2},$
$H_{(5;3)}^{12}=N_{(5;3)}^{3},$  $H_{(5;3)}^{13}=N_{(5;3)}^{4},$
$H_{(9;3)}^{14}=K_{(9;3)}. $

$(iii) $ Also, applying the program $ BGroidAP3 $ for $ K_{(9;3)}
,$ we obtain that it has $ 5 $ wide subgroupoids. We observe that,
each wide subgroupoid of this groupoid is a normal subgroupoid.
\hfill$.\b$

{\it  Example 5.2.} $(i)$ {\it Determination of normal subgroupoids
of the groupoid $ G_{(8;2)}. $} Use the inputs data presented in the
Example 3.2 and execute the program $ BGroidAP4. $ This groupoid has
the following $ 5 $ normal subgroupoids:
$N_{(2;2)}^{1}=G_{(8;2),0},$  $N_{(4;2)}^{2}=\{ g_{1}, g_{2}, g_{3},
g_{8}\},$ $N_{(4;2)}^{3}=\{ g_{1}, g_{2}, g_{4}, g_{6} \},$
$N_{(4;2)}^{4}=\{ g_{1}, g_{2}, g_{5}, g_{7}\},$ $N_{(8;2)}^{7}=
G_{(8;2)}.$

$(ii) $ {\it Determination of normal subgroupoids of the groupoid
$ K_{4}\coprod {\cal F}_{(4;2)}({\bf R}^{2}). $} Applying the
program $ BGroidAP4 $ for $ G = K_{4}\coprod {\cal F}_{(4;2)}({\bf
R}^{2}), $ we obtain that $ G $ has $ 10 $ normal subgroupoids,
namely: $N_{(3;3)}^{1}=G_{0},$  $N_{(4;3)}^{2}=\{(1),Id_{Ox},
Id_{Oy},$ $\sigma\},$ $N_{(4;3)}^{3}=\{(1),Id_{0x},
Id_{Oy},\tau\},$ $N_{(4;3)}^{4}=\{(1),Id_{Ox},
Id_{Oy},\sigma\circ\tau\},$ $N_{(5;3)}^{5}=\{ (1), Id_{Ox},
Id_{0y}, \sigma_{Ox}, \sigma_{Oy}\},$ $N_{(6;3)}^{6}=\{ (1),
Id_{Ox},Id_{Oy},\sigma,\tau,\sigma\circ \tau\},$ $N_{(6;3)}^{7}=\{
(1), Id_{Ox},Id_{Oy},\sigma,\sigma_{Ox},\sigma_{Oy}\},$
$N_{(6;3)}^{8}=\{ (1),
Id_{Ox},Id_{Oy},\tau,\sigma_{Ox},\sigma_{Oy}\},$ $N_{(6;3)}^{9}=\{
(1), Id_{Ox},Id_{Oy},\sigma\circ \tau,\sigma_{Ox},\sigma_{Oy}\},$
$N_{(8;3)}^{10}= G.$\hfill$\b$

{\it Example 5.3.} {\it Determination of subgroups and normal
subgroups of a finite group.} We consider the dihedral group $
D_{5}=\{ x_{1} = e, x_{2} = a, x_{3} = a^{2}, x_{4} = a^{3}, x_{5} =
a^{4}, x_{6} =b, x_{7} = ab, x_{8} = a^{2}b, x_{9} = a^{3}b, x_{10}
= a^{4}b\} $
 generated by the elements $ a, b $ with properties
$ a^{5}= e $ and $ b^{2} = e. $

The inputs data for this group are the following:\\[-0.4cm]
$$\begin{array}{cccccccccc} \\
10 & & & & & & & & &\\
1 & & & & & & & & &\\[0.1cm]

1 & 1 & 1 & 1 & 1 & 1 & 1 & 1 & 1 & 1\\
1 & 1 & 1 & 1 & 1 & 1 & 1 & 1 & 1 & 1\\
1 & 5 & 4 & 3 & 2 & 6 & 7 & 8 & 9 & 10\\[0.1cm]
\end{array}   \begin{array}{cccccccccc} \\
 1 & 2 & 3 & 4 & 5 & 6 & 7 & 8 & 9 & 10\\
 2 & 3 & 4 & 5 & 1 & 7 & 8 & 9 & 10 & 6\\
 3 & 4 & 5 & 1 & 2 & 8 & 9 & 10 & 6 & 7\\
 4 & 5 & 1 & 2 & 3 & 9 & 10 & 6 & 7 & 8\\
 5 & 1 & 2 & 3 & 4 & 10 & 6 & 7 & 8 & 9\\
 6 & 10 & 9 & 8 & 7 & 1 & 5 & 4 & 3 & 2\\
 7 & 6 & 10 & 9 & 8 & 2 & 1 & 5 & 4 & 3\\
 8 & 7 & 6 & 10 & 9 & 3 & 2 & 1 & 5 & 4\\
 9 & 8 & 7 & 6 & 10 & 4 & 3 & 2 & 1 & 5\\
10 & 9 & 8 & 7 & 6 & 5 & 4 & 3 & 2  & 1\\
\end{array}$$\\[-0.5cm]

Use the above input data and execute the program $ BGroidAP4. $
Then $ D_{5} $ has $ 3 $ normal subgroups, namely:$ N_{1} =\{e\},
N_{2}=\{ e, a, a^{2}, a^{3}, a^{4}\} $ and $ N_{3}=D_{5}. $

Applying the program $ BGroidAP2 $ or $ BGroidAP3 ,$ we obtain
that $ D_{5} $ has $ 8 $ subgroups ( in fact, subgroupoids and
wide subgroupoids with one unit ), namely: $ H_{1}=\{e\},
H_{2}=\{e,b\}, H_{3}=\{e,ab\}, H_{4}=\{e,a^{2}b\},
H_{5}=\{e,a^{3}b\}, H_{6}=\{e,a^{4}b\}, H_{7}=\{e, a, a^{2},
a^{3}, a^{4}\}, H_{8}=D_{5}$\hfill$.\b$

For more details concerning the programs given in this paper, the reader can be inform at
e-mail adress: ivan@math.uvt.ro.\\[-0.4cm]

\begin{center}
{\bf REFERENCES}
\end{center}

1. R. Brown, From Groups to Groupoids : a brief survey, {\it Bull.
London Math. Soc.},{\bf 19}(1987), 113-134.

2. A. Coste, P. Dazord and A. Weinstein, Groupoides symplectiques,
 {\it Publ. Dept. Math. Lyon,} 2/A (1987),1-62.

3. Gh. Ivan, Algebraic constructions of Brandt groupoids, {\it
Proceedings of the Algebra Symposium, "Babe\c s-Bolyai" University,
Cluj,} (2002), 69-90.

4. Gh. Ivan and G. Stoianov, The program $ BGroidAP2 $ for
determination of all subgroupoids of a Brandt groupoid, {\it Analele
Univ. de Vest, Timi\c soara, Seria Matematic\u a- Informatic\u a,}
vol. {\bf 42}, fasc. 1 (2004), 93-117.

5. M. Ivan,  General properties of the symmetric groupoid of a
finite set, {\it An. Univ. Craiova, Ser. Mat. Informatic\u a,} {\bf
30} (2003),no.2, 109-119.

6. P. J. Higgins, {\it Notes on Categories and Groupoids}, Von
Nostrand Reinhold, London, 1971.

7. K. Mackenzie, {\it Lie Groupoids and Lie Algebroids in
Differential Geometry}, London Math. Soc., Lect. Notes Series, {\bf
124}, Cambridge Univ.Press, 1987.

8. A. Weinstein, Groupoids: Unifying Internal and External
Symmetries, {\it Notices Amer. Math. Soc.,} {\bf 43} (1996), 744-752.\\

\end{document}